\documentclass[final,5p,times,twocolumn]{elsarticle}

\usepackage{amssymb}
\usepackage{amsmath}
\usepackage{color}

\newcommand{\R}{\mathbb{R}}

\newcommand{\N}{\mathbb{N}}

\newtheorem{thm}{Theorem}
\newtheorem{rem}{Remark}
\newtheorem{prop}{Proposition}
\newtheorem{cor}{Corollary}
\newcommand{\change}[1]{#1}

\journal{Systems \& Control Letters}

\begin{document}

\begin{frontmatter}


\title{Improved residual mode separation for finite-dimensional control of PDEs:\\ Application to the Euler--Bernoulli beam}

\author[label1]{Anton Selivanov}
\affiliation[label1]{organization={School of Electrical and Electronic Engineering, The University of Sheffield},
             city={Sheffield},
             country={UK}}

\author[label2]{Emilia Fridman} 

\affiliation[label2]{organization={School of Electrical Engineering, Tel Aviv University},
            city={Tel Aviv},
            country={Israel}}

\begin{abstract}
We consider a simply-supported Euler--Bernoulli beam with viscous and Kelvin--Voigt damping. Our objective is to attenuate the effect of an unknown distributed disturbance using one piezoelectric actuator. We show how to design a {state-feedback controller based on a finite number of dominating modes that guarantees that the $L^2$ gain is not greater than a given value}. If the remaining (infinitely many) modes are simply ignored, the calculated $L^2$ gain is wrong. This happens because of the spillover phenomenon that occurs when the effect of the control on truncated modes is not accounted for in the feedback design. We propose a simple modification of the $H_\infty$ cost that prevents spillover. The key idea is to treat the control as a disturbance in the truncated modes and find the corresponding $L^2$ gains using the bounded real lemma. These $L^2$ gains are added to the control weight in the $H_\infty$ cost for the dominating modes, which prevents spillover. A numerical simulation of an aluminum beam with realistic parameters demonstrates the effectiveness of the proposed method. The presented approach is applicable to other types of PDEs, such as the heat, wave, and Kuramoto–Sivashinsky equations, as well as their semilinear versions. {The proposed method gives a Lyapunov functional that can also be used for guaranteed cost control, regional stability analysis, and input-to-state stability.}
\end{abstract}


\end{frontmatter}



\section{Introduction}
The $H_\infty$ control theory enables the design of controllers robust to modeling errors, measurement noise, and unknown disturbances \cite{Francis1987,Green2012}. Its extension to \textit{infinite}-dimensional systems is challenging, especially if the controller is required to be \textit{finite}-dimensional. In particular, the direct extension of the \textit{frequency-domain} approach results in infinite-dimensional controllers \cite{ozbay1990,ozbay1993a}, which are then approximated by finite-dimensional ones~\cite{Nett1983,Glover1986,ozbay1991}. The approximation leads to performance degradation known as spillover, which can be characterized via the $H_\infty$ norm of the approximation error \cite{Bontsema1988}. The frequency-domain approach, which is very natural for $H_\infty$ control, may be difficult to use in the infinite-dimensional case since the transcendental transfer function of an infinite-dimensional plant may be hard to find \cite{Curtain1995}, its inner-outer factorization required for the design is not straightforward \cite{ozbay1993a}, and the MIMO case requires restrictive assumptions \cite{ozbay1990}. Furthermore, the frequency-domain approach is not applicable to guaranteed cost control, regional stability analysis, systems with time-varying delays, and nonlinear systems. 


The \textit{time-domain} approach avoids these restrictions. Its direct extension to infinite-dimensional systems leads to an operator Riccati equation, which also results in an infinite-dimensional controller~\cite{Curtain1990a,Keulen1993}. To obtain a finite-dimensional controller, one can perform modal decomposition~\cite{Athans1970,Balas1982a} and design a controller for a finite number of dominating modes~\cite{Triggiani1980,Balas1982,Curtain1993}. Similarly to the frequency-domain design, this leads to spillover: neglected modes deteriorate the overall system performance~\cite{Balas1978a,Meirovitch1983}. Nevertheless, \textit{stability} under spillover can be guaranteed using residual filters \cite{Balas1988,Moheimani1998,Harkort2011} or spectral properties of linear operators representing the dynamics \cite{Lasiecka1983,Curtain1984}. 

The time-domain \textit{performance} analysis under spillover is more challenging than the stability analysis and requires careful treatment of the neglected modes. Such treatment has been provided for parabolic PDEs in \cite{Karafyllis2016a,Karafyllis2019k,Katz2020a,Karafyllis2021} with subsequent extensions to input/output delays \cite{Selivanov2018e,Lhachemi2021,Katz2022}, semilinear systems \cite{Katz2021c,Selivanov2024a}, as well as the Kuramoto--Sivashinsky~\cite{Katz2022a}, wave \cite{Selivanov2023b}, and Euler--Bernoulli \cite{Selivanov2023} equations. 

This paper proposes a new way of dealing with spillover in the time domain. Namely, we treat the control input as a disturbance in the residue modes and explicitly solve the algebraic Riccati equation for each neglected mode to find the input-to-state $L^2$ gains. These gains are added to the control weight in the cost used to design {a state-feedback controller that guarantees that the $L^2$ gain is not greater than a given value.} This idea leads to a simple yet efficient way of designing a finite-dimensional controller that avoids spillover. The analysis is based on the cost decomposition presented in Section~\ref{sec:cost decomposition}. We develop this idea to attenuate disturbances in the Euler--Bernoulli beam with piezoelectric actuators, which is of great importance for aerospace, civil, and mechanical engineering. Using a numerical example of an aluminum beam, we demonstrate a drastic improvement compared to our previous results in~\cite{Selivanov2023}. Namely, spillover is avoided using just $8$ modes instead of $32$, and we prove that the $L^2$ gain can only decrease when more modes are considered. 

The frequency-domain approach to the $H_\infty$ control of beams was developed in \cite{Bontsema1988a,Lenz1993,Halim2001}, the controllability problem under piezoelectric actuators was studied in \cite{Tucsnak1996,Crepeau2006,Bai2024}, and experimental results (without spillover analysis) were reported in \cite{Halim2002,Belyaev2018}. Here, we develop the time-domain method, which, differently from the frequency-domain approach, can be extended to {guaranteed cost control and regional stability analysis}. Furthermore, the proposed idea can significantly improve the finite-dimensional controller design for other types of PDEs, including the heat, wave, and Kuramoto--Sivashinky equations. It also admits an extension to semilinear PDEs in a manner similar to \cite{Katz2021c,Selivanov2024a}.

\textit{Notations:} $|\cdot|$ is the Euclidean norm, $\|\cdot\|$ is the $L^2$ norm, $\langle\cdot,\cdot\rangle$ is the scalar product in $L^2$, $H^p(0,\pi)$ with $p\in\N$ are the Sobolev spaces, $H^{-p}(0,\pi)$ are their dual spaces, {$H_{loc}^1(0,\infty)$ are the functions that belong to $H^1(K)$ for any compact $K\subset(0,\infty)$,} $\operatorname{diag}\{\omega_1,\ldots,\omega_N\}$ is the diagonal matrix with diagonal elements $\omega_n$, $n=1,\ldots,N$. For a matrix $P$, the notation $P<0$ implies that $P$ is square, symmetric, and negative-definite. Partial derivatives are denoted by indices, e.g., $z_t=\partial z/\partial t$. 
\subsection{Preliminaries: $H_\infty$ control of finite-dimensional systems} 
Consider the LTI system 
\begin{equation}\label{LTI}
\begin{aligned}
    \dot x(t)&=Ax(t)+Bu(t)+Ev(t),\qquad x(0)=0,\\
    y(t)&=Cx(t)+Du(t)
\end{aligned}
\end{equation}
with state $x\in\R^n$, control input $u\in\R^m$, disturbance $v\in\R^k$, controlled output $y\in\R^l$, and constant matrices $A$, $B$, $C$, $D$, and $E$. {We say that $u=-Kx$ with $K\in\R^{m\times n}$ guarantees that the $L^2$ gain from $v$ to $y$ is not greater than $\gamma>0$ if the solutions of the closed-loop system satisfy}
\begin{multline}\label{HinfProblem}
\textstyle\int_0^\infty\left[|y(t)|^2-\gamma^2|v(t)|^2\right]dt\le0,\ 
\forall v\in L^2([0,\infty),\R^k).
\end{multline}
{That is, the $H_\infty$ norm of the closed-loop transfer function from $v$ to $y$ is not greater than $\gamma$.} The proofs of the following results are given, e.g., in \cite{Green2012}.
\begin{prop}\label{prop:HinfDesign}
Consider \eqref{LTI} such that $D^\top C=0$ and $R=D^\top D>0$. Given $\gamma>0$, let $0<P\in\R^{n\times n}$ satisfy 
\begin{equation}\label{ARE}
   PA+A^\top P-P(BR^{-1}B^\top-\gamma^{-2}EE^\top)P+C^\top C=0. 
\end{equation}
Then $u(t)=-R^{-1}B^\top Px(t)$ guarantees \eqref{HinfProblem}. 
\end{prop}
\begin{rem}[Solution existence]\label{rem:AREexistence}
If $(A,B)$ is stabilizable, $(A,C)$ is detectable, and $\gamma$ is large enough, then \eqref{ARE} has a solution. For this solution, the closed-loop matrix $A-BR^{-1}B^\top P$ is stable. 
\end{rem}
\begin{cor}[Bounded Real Lemma]\label{cor:BRL}
Consider \eqref{LTI} with $B=0$ and $D=0$ (i.e., without control). Given $\gamma>0$, let $0<P\in\R^{n\times n}$ satisfy 
\begin{equation}\label{BRL:ARE}
    PA+A^\top P+\gamma^{-2}PEE^\top P+C^\top C=0. 
\end{equation}
Then \eqref{HinfProblem} holds without control. 
\end{cor}
\section{Model description}
\subsection{Euler--Bernoulli beam with control and disturbance}
We consider the Euler--Bernoulli beam described by 
\begin{equation}\label{EB_beam_original}
\begin{aligned}
&\mu\tilde z_{tt}(x,t)+EI\tilde z_{xxxx}(x,t)+c_v\tilde z_t(x,t)+c_kI\tilde z_{xxxxt}(x,t)=\\
&\hspace{2.5cm}c_a[\delta^{\prime}(x-\tilde x_L)-\delta^{\prime}(x-\tilde x_R)]\tilde u(t)+\tilde w(x,t),\\
&\tilde z(0,t)=\tilde z_{xx}(0,t)=\tilde z(L,t)=\tilde z_{xx}(L,t)=0, 
\end{aligned}
\end{equation}
where $\tilde z\colon[0,L]\times[0,\infty)\to\R$ is the transverse deflection of a beam of length $L$, linear density $\mu$, Young's modulus of elasticity $E$, and moment of inertia $I$. The model accounts for the viscous damping $c_v\tilde z_t$ and structural (Kelvin--Voigt) damping $c_kI\tilde z_{xxxxt}$ \cite{Russell1992,Herrmann2008}. The external disturbance is represented by $\tilde w\colon(0,L)\times[0,\infty)\to\R$. All the parameters are constant in time and space. The boundary conditions correspond to the hinged ends. 

A piezoelectric actuator produces bending moment on $[\tilde x_L,\tilde x_R]\subset(0,L)$ proportional to the applied voltage $\tilde u\colon[0,\infty)\to\R$. Namely, $m(x,t)=c_a[h(x-\tilde x_L)-h(x-\tilde x_R)]\tilde u(t)$, where $h(x)$ is the step function. Since $m(x,t)$ contributes to the beam's moment-curvature relationship, it enters \eqref{EB_beam_original} through the Laplace operator: $m_{xx}(x,t)=c_a[\delta'(x-\tilde x_L)-\delta'(x-\tilde x_R)]\tilde u(t)$, where $\delta'(\,\cdot\,-\tilde x)\in H^{-2}(0,L)$ is the derivative of the Dirac delta function defined as 
\begin{equation}\label{deltaDefinition}
    \delta'(\cdot-\tilde x)f=\int_0^L\delta'(x-\tilde x)f(x)\,dx=-f'(\tilde x) 
\end{equation}
for any $\tilde x\in(0,L)$ and $f\in H^2(0,L)$. Note that it is natural for the derivatives of $\delta$ to be of opposite signs since the piezoelectric patch applies forces of opposite directions to its ends when it contracts or expands. A more detailed study of piezoelectric actuators is provided in \cite{Crawley1987,Fuller1996,Halim2001}. 
\begin{rem}[Damping model]\label{rem:squareRoot}
The Kelvin--Voigt damping is motivated by the experimental observation that damping rates in beams increase with frequency \cite{Russell1992}. This is also captured by the ``square root'' model given by $-c_r\tilde z_{xxt}$ \cite{Chen1982}. Our analysis can be extended to the ``square root'' model straightforwardly. 
\end{rem}
By scaling the space and time as follows
\begin{equation}\label{changeOfVariables}
    z(x,t)=\tilde z(a_1x,a_2t),\quad 
    a_1=\frac{L}{\pi},\quad 
    a_2=a_1^2\sqrt{\frac{\mu}{EI}},
\end{equation}
we rewrite \eqref{EB_beam_original} as 
\begin{equation}\label{EB_beam}
\begin{aligned}
& z_{tt}+z_{xxxx}+c_1 z_t+c_2 z_{xxxxt}=\left[\delta_L'-\delta_R'\right]u+w, \\
& \change{z(0,t)=z_{xx}(0,t)=z(\pi, t)=z_{xx}(\pi,t)=0,}
\end{aligned}
\end{equation}
where $x\in[0,\pi]$, $t\ge0$, 
\begin{equation*}
\begin{aligned}
&c_1=\frac{c_va_2}{\mu},\quad c_2=\frac{c_kIa_2}{\mu a_1^4},\quad x_L=\frac{\tilde x_L}{a_1},\quad x_R=\frac{\tilde x_R}{a_1}, \\
&\delta_L'=\delta'\left(x-x_L\right),\quad\delta_R'=\delta'\left(x-x_R\right),\\
&u(t)=\frac{c_aa_2^2}{\mu a_1^2}\tilde u(a_2t),\quad w(x,t)=\frac{a_2^2}{\mu}\tilde w(a_1x,a_2t). 
\end{aligned}
\end{equation*}
Note that \eqref{deltaDefinition} implies $\delta'(a_1x-\tilde x)=\delta'(x-\tilde x/a_1)/a_1^2$. To simplify further derivations, we assume that 
\begin{equation}\label{c1c2cond}
    c_1+c_2\le\sqrt{2}. 
\end{equation}
That is, the dynamics are dominated by the elasticity rather than damping. The extension to $c_1+c_2>\sqrt{2}$ is straightforward. 
\subsection{Well-posedness}\label{sec:Well-posedness}
\change{Let $H_{BC}^2$ and $H_{BC}^4$ be the closure in $H^2$ and $H^4$, respectively, of all functions $f\in C^\infty[0,\pi]$ satisfying $f^{(2k)}(0)=0=f^{(2k)}(\pi)$ for all $k\ge0$.}
The energy space of \eqref{EB_beam} is 
\begin{equation*}
X=H_{BC}^2(0,\pi)\times L^2(0,\pi)
\end{equation*}
with the scalar product 
\begin{equation*}
\langle (f_1,g_1),(f_2,g_2)\rangle_X=\langle f_1'',f_2''\rangle_{L^2}+\langle g_1,g_2\rangle_{L^2}. 
\end{equation*}
Consider  
\begin{equation}\label{A0}
\mathcal{A}_0f=-f'',\quad D(\mathcal{A}_0)=H^2_{BC}(0,\pi)\subset L^2(0,\pi). 
\end{equation}
In the operator form, \eqref{EB_beam} is written as 
\begin{equation}\label{abstractPDE}
\dot{\bar z}=\mathcal{A}\bar z+f, 
\end{equation}
where 
\begin{equation*}
\begin{aligned}
&\bar z(t)=\begin{bmatrix}
z(\cdot, t)\\z_t(\cdot, t)
\end{bmatrix},\quad
\mathcal{A}=\begin{bmatrix}
0 & I \\ -\mathcal{A}_0^2 & -(c_1I+c_2\mathcal{A}_0^2)
\end{bmatrix},\\
&f(t)=\begin{bmatrix}
0\\ [\delta_L'-\delta_R']u(t)+w(\cdot,t)
\end{bmatrix}. 
\end{aligned}
\end{equation*}
Since $D(\mathcal{A}_0^2)=H^4_{BC}(0,\pi)$, we have 
\begin{equation*}
D(\mathcal{A})=X_1=H_{BC}^4(0,\pi)\times H_{BC}^4(0,\pi)\subset X. 
\end{equation*}
The adjoint of $\mathcal{A}$ with respect to the scalar product in $X$ is 
\begin{equation*}
\mathcal{A}^*=\begin{bmatrix}
0 & -I \\ \mathcal{A}_0^2 & -(c_1I+c_2\mathcal{A}_0^2)
\end{bmatrix},\quad D(\mathcal{A}^*)=X_1\subset X. 
\end{equation*}
{
Both $\mathcal{A}$ and $\mathcal{A}^*$ are dissipative. To see this, consider 
\begin{multline*}
\langle \mathcal{A}(f,g),(f,g)\rangle_X=\langle (g,-f^{(4)}-c_1g-c_2g^{(4)}),(f,g)\rangle_X\\
=\langle g'',f''\rangle_{L^2}-\langle f^{(4)},g\rangle_{L^2}-c_1\langle g,g\rangle_{L^2}-c_2\langle g^{(4)},g\rangle_{L^2}. 
\end{multline*}
For $f,g\in H_{BC}^4$, integration by parts gives 
\begin{equation*}
    \langle f^{(4)},g\rangle_{L^2}=\langle f'',g''\rangle_{L^2}\quad\text{and}\quad
    \langle g^{(4)},g\rangle_{L^2}=\langle g'',g''\rangle_{L^2}. 
\end{equation*}
Substituting, we obtain 
\begin{equation*}
    \langle \mathcal{A}(f,g),(f,g)\rangle_X=-c_1\|g\|_{L^2}^2-c_2\|g''\|_{L^2}^2\le0. 
\end{equation*}
Similarly, 
\begin{equation*}
    \langle \mathcal{A}^*(f,g),(f,g)\rangle_X=-c_1\|g\|_{L^2}^2-c_2\|g''\|_{L^2}^2\le0. 
\end{equation*}}%
Since {$\overline{D(\mathcal{A})}=X$, $\mathcal{A}$ is closed, and }$\mathcal{A}$ and $\mathcal{A}^*$ are dissipative, $\mathcal{A}$ generates a $C_0$-semigroup of contractions on $X$ \cite[Corollary~4.4]{Pazy1983}.  

The set $D(\mathcal{A}^*)$ with the norm ${\|z\|_D=\|(\bar\beta I-\mathcal{A}^*)z\|_X}$, where $\beta$ is any regular point of $\mathcal{A}$ and $\bar\beta$ is its complex conjugate, is a Hilbert space \cite[Proposition 2.10.1]{Tucsnak2009}. \change{Its dual with respect to the pivot space $X$ is $X_{-1}=Y\times \change{H_{BC}^{-4}(0,\pi)}$, where $Y$ and $H_{BC}^{-4}$ are the dual spaces of $H_{BC}^4$ with respect to $H^2$ and $L^2$, respectively.}

We assume that 
\begin{equation}\label{w condition}
w\in H^1_{loc}((0,\infty),\change{H_{BC}^{-4}(0,\pi)})\cap L^2((0,\infty),L^2(0,\pi)). 
\end{equation}
The control input that we design later satisfies $u\in H^1_{loc}((0,\infty),\R)$. \change{Since $H_{BC}^4\subset H^2$, we have $\delta_L',\delta_R'\in H_{BC}^{-4}\supset H^{-2}$. Since $\delta_L'$ and $\delta_R'$ are constant in time, this implies }
\begin{equation*}
[\delta_L'-\delta_R']u\in H^1_{loc}((0,\infty),\change{H^{-4}_{BC}(0,\pi)}). 
\end{equation*}
Therefore, $f\in H^1_{loc}((0,\infty),X_{-1})$. By \cite[Theorem~4.1.6]{Tucsnak2009}, for $z(\cdot,0)\in H^2_{BC}(0,\pi)$ and $z_t(\cdot,0)\in L^2(0,\pi)$, there exists a unique solution of \eqref{abstractPDE} in $X_{-1}$ that satisfies 
\begin{equation*}
    \bar z\in C([0,\infty),X)\cap C^1([0,\infty),X_{-1}). 
\end{equation*}
{Since $\bar z=(z,z_t)^T$, this implies 
\begin{equation*}
z\in C([0,\infty),H_{BC}^2(0,\pi))\quad\text{and}\quad 
z_t\in C([0,\infty),L^2(0,\pi)). 
\end{equation*}}%
\section{{Robust state-feedback} control of the beam}
Given non-negative scalars $\rho_x$, $\rho_u$, and $\gamma$, our objective is to find a state-feedback control law guaranteeing that the trajectories of \eqref{EB_beam} with $z(\cdot,0)\equiv0\equiv z_t(\cdot,0)$ satisfy (cf.~\eqref{HinfProblem})
\begin{multline}\label{HinfObjective}
    J=\int_0^\infty\Bigl[\|z(\cdot,t)\|^2+\rho_x\|z_{xx}(\cdot,t)\|^2\\
        +\rho_uu^2(t)-\gamma^2\|w(\cdot,t)\|^2\Bigr]dt\le0
\end{multline}
for all $w$ satisfying \eqref{w condition}. Such control guarantees that the $L^2$ gain is not greater than $\gamma$. {Using \eqref{changeOfVariables}, one can rewrite \eqref{HinfObjective} in terms of the original state, input, and disturbance with $\tilde\rho_x=\rho_xa_1^4$, $\tilde\rho_u=\rho_u\frac{c_a^2a_2^4}{\mu^2a_1^3}$, and $\tilde\gamma=\gamma\frac{a_2^2}{\mu}$.}
\begin{rem}[Performance index]
Since the potential energy of \eqref{EB_beam_original} due to bending is $\frac{EI}2\|\tilde z_{xx}(\cdot,t)\|^2$ \cite[p.~317]{Timoshenko1955}, we include $\|z_{xx}(\cdot,t)\|^2$ in \eqref{HinfObjective}. The kinetic energy of \eqref{EB_beam_original} is $\frac{\mu}{2}\|\tilde z_t(\cdot,t)\|^2$, and it is natural to include $\|z_t(\cdot,t)\|^2$ in \eqref{HinfObjective}. To simplify the exposition, we do not present this extension, which requires one to consider multiple cases depending on the values of $c_1$, $c_2$, $\rho_x$, and $\rho_u$. 
\end{rem}
\subsection{Modal decomposition}\label{sec:modal decomposition}
The modes and natural frequencies of \eqref{EB_beam} are 
\begin{equation*}
    \varphi_n(x)=\sqrt{2/\pi}\sin nx,\quad \omega_n=n^2,\qquad n\in\N. 
\end{equation*}
These are the eigenfunctions and eigenvalues of $\mathcal{A}_0$ defined in \eqref{A0}, which form a complete orthonormal system in $L^2(0,\pi)$. Therefore, 
\begin{equation*}
    z(\cdot,t)\stackrel{L^2}{=}\sum_{n=1}^\infty z_n(t)\varphi_n,\qquad z_n(t)=\langle z(\cdot,t),\varphi_n\rangle. 
\end{equation*}
Substituting this into \eqref{EB_beam}, in view of 
\begin{equation*}
	\langle\delta_L',\varphi_n\rangle\stackrel{\eqref{deltaDefinition}}{=}-\varphi_n'(x_L)\quad\text{and}\quad
	\langle\delta_R',\varphi_n\rangle\stackrel{\eqref{deltaDefinition}}{=}-\varphi_n'(x_R), 
\end{equation*}
we obtain the ODEs for the Fourier coefficients
\begin{equation*}
    \ddot{z}_n(t)+2\zeta_n\omega_n\dot z_n(t)+\omega_n^2z_n(t)=b_nu(t)+w_n(t),\quad n\in\N, 
\end{equation*}
where 
\begin{equation*}
\begin{aligned}
&\zeta_n=(c_1\omega_n^{-1}+c_2\omega_n)/2,\\
&b_n=n\sqrt{2/\pi}\left(\cos n x_R-\cos n x_L\right),\\
&w_n(t)=\left\langle w(\cdot,t), \varphi_n\right\rangle. 
\end{aligned}
\end{equation*}
The ODEs can be written as 
\begin{equation}\label{FourierCoefficients}
    \dot{\bar z}_n(t)=A_n \bar{z}_n(t)+B_nu(t)+E_nw_n(t), \quad n \in \mathbb{N}, 
\end{equation}
where 
\begin{equation*}
\bar z_n=\begin{bmatrix}
z_n\\\dot z_n
\end{bmatrix},\ 
A_n=\begin{bmatrix}
0 & 1 \\ -\omega\ _n^2 & -2\zeta_n\omega_n
\end{bmatrix},
B_n=\begin{bmatrix}
0 \\  b_n
\end{bmatrix},\ 
E_n=\begin{bmatrix} 0 \\ 1 \end{bmatrix}. 
\end{equation*}
\begin{figure}
\includegraphics[width=\linewidth]{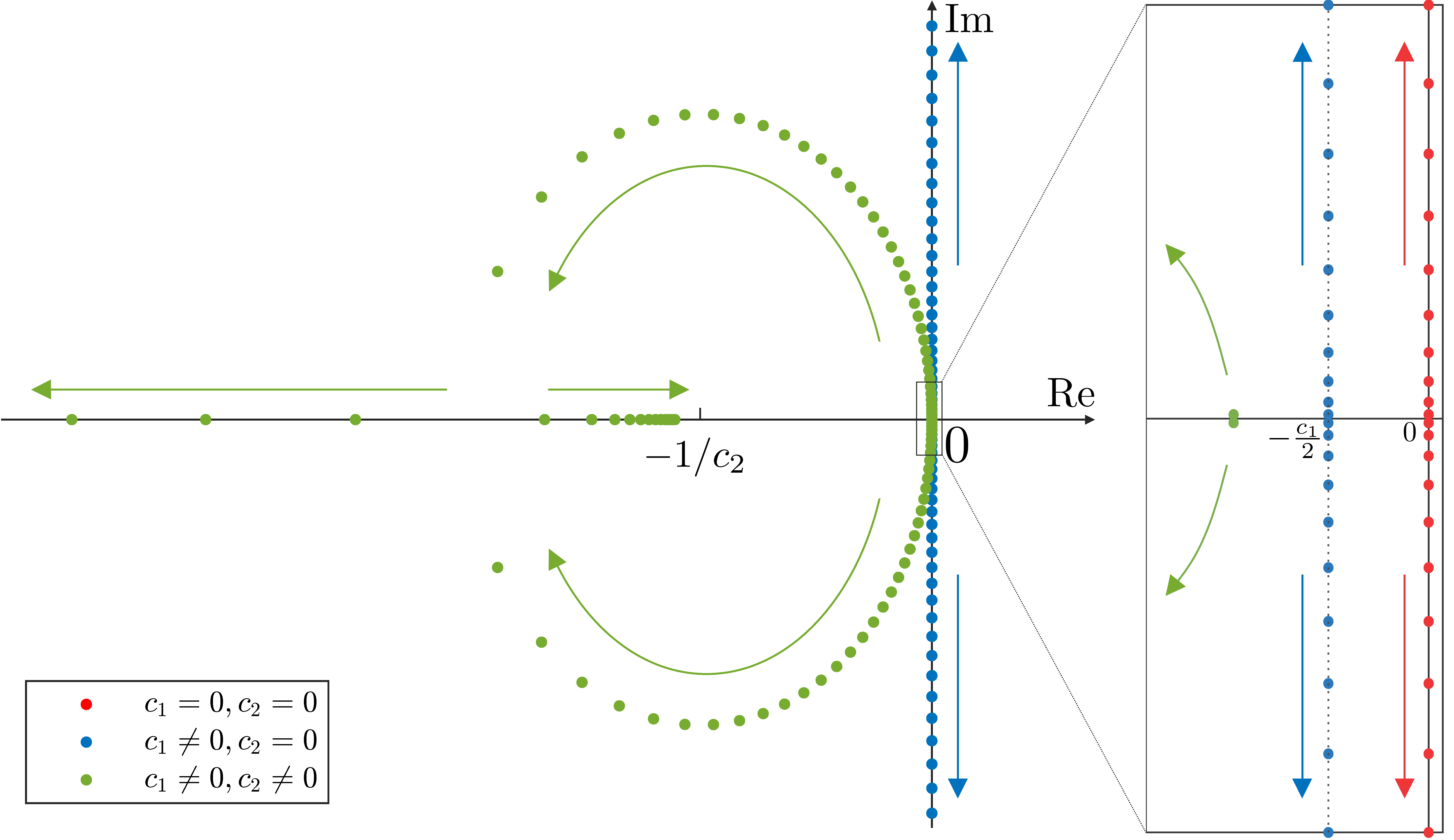}
\caption{The eigenvalues of $A_n$ given in \eqref{eigenvalues} for $n=1,\ldots,50$. Red dots --- no damping ($c_1=0=c_2$); blue dots --- viscous damping ($c_1=1.4\times 10^{-3}$, $c_2=0$); green dots --- viscous and Kelvin--Voigt damping ($c_1=1.4\times 10^{-3}$, $c_2=1.3\times 10^{-3}$).}\label{fig:eigenvalues}
\end{figure}%
The eigenvalues of $A_n$ are 
\begin{equation}\label{eigenvalues}
    \lambda_n^{\mp}=-\omega_n(\zeta_n\pm\sqrt{\zeta_n^2-1}). 
\end{equation}
Without damping ($\zeta_n=0$), infinitely many imaginary roots $\lambda_n^{\pm}=\pm i\omega_n$ (see Fig.~\ref{fig:eigenvalues}) give rise to free vibrations in \eqref{EB_beam} in the absence of control and disturbance. Our approach does not work in this case since it is not enough to deal only with a finite number of modes. Viscous damping ($c_1\neq 0$) ensures that $\operatorname{Re}\lambda_n^\pm=-c_1/2$. Kelvin--Voigt damping ($c_2\neq 0$) improves the stability further guaranteeing  
\begin{equation*}
\operatorname{Re}\lambda_n^-\to-\infty\quad\text{and}\quad\operatorname{Re}\lambda_n^+\to-1/c_2. 
\end{equation*}
We develop our approach for the case when $c_1\neq0\neq c_2$. 
\subsection{The spillover phenomenon}\label{sec:spillover}
It is common in engineering practice to design controllers based on a few dominating modes while ignoring the residue. This subsection demonstrates that such an approach may suffer from the spillover phenomenon.

\begin{figure}
\includegraphics[width=\linewidth]{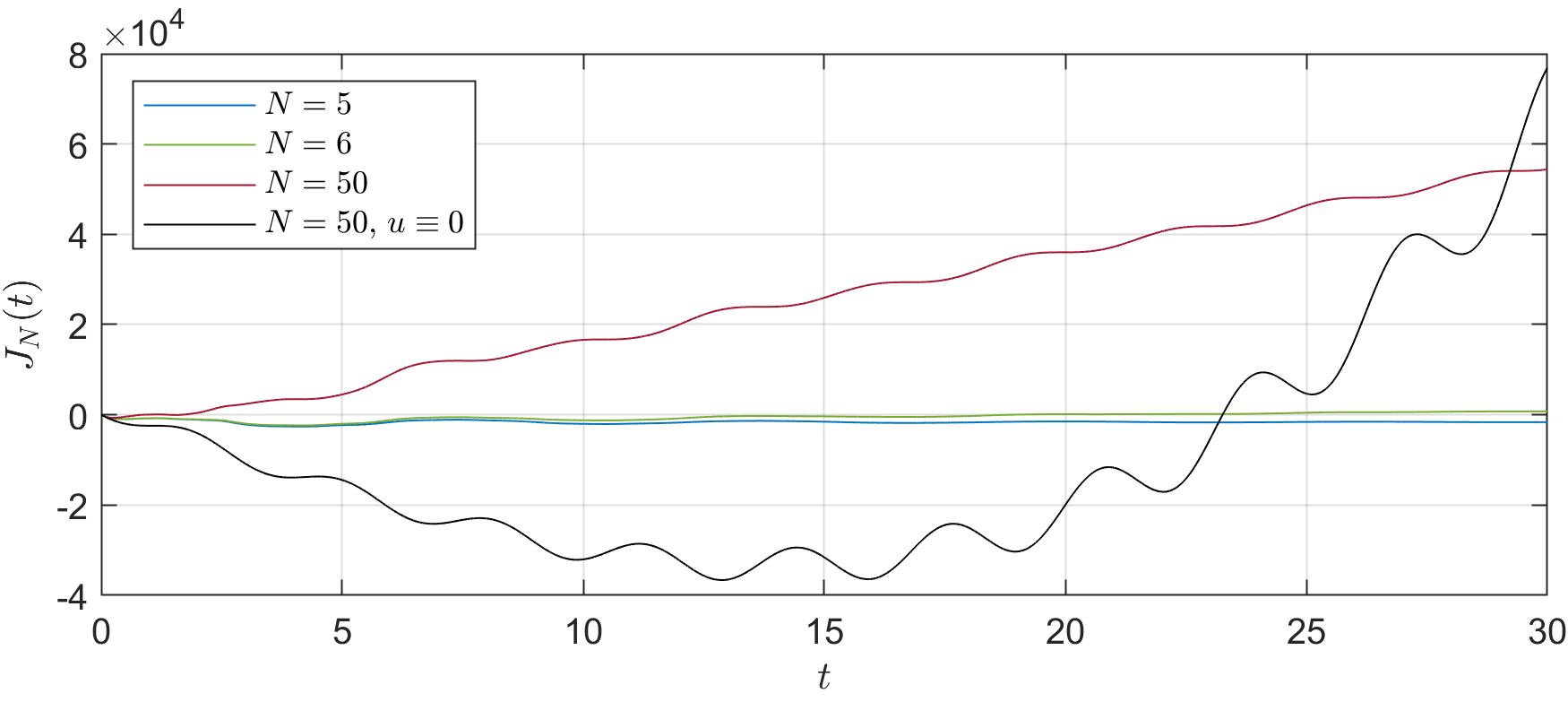}
\caption{The value of $J_N(t)$, defined in \eqref{J(t)}, for $N=5$ (blue), $N=6$ (green), and $N=50$ (red). The black line is $J_N(t)$ for $N=50$ and $u\equiv 0$. A controller designed for the first $5$ modes cannot guarantee \eqref{HinfObjective} for the original system because of the spillover phenomenon.}\label{fig:spillover}
\end{figure}%
Consider the Euler--Bernoulli beam \eqref{EB_beam} with 
\begin{equation*}
	c_1=1.4\times 10^{-3},\ c_2=1.3\times 10^{-3},\ x_L=0.91,\ x_R=0.97. 
\end{equation*}
The choice of the parameters is explained in Section~\ref{sec:example}. Let us try to design a controller guaranteeing \eqref{HinfObjective} with $\rho_x=0.1$ and $\rho_u=10^{-3}$ by considering only $5$ modes in the modal decomposition \eqref{FourierCoefficients}. Using Proposition~\ref{prop:HinfDesign} (see Section~\ref{sec:HinfDesign} for details), we find $\gamma\approx6.97$ and the corresponding controller gain $-R^{-1}B^\top P\in\R^{1\times 10}$. 

Figure~\ref{fig:spillover} shows the values of 
\begin{multline}\label{J(t)}
    \textstyle J_N(t)=\int_0^t\Bigl[\sum_{n=1}^N\left[(1+\rho_x\omega_n^2)z_n^2(t)-\gamma^2w_n^2(t)\right]\\
    \textstyle+\rho_uu^2(t)\Bigr]dt 
\end{multline}
for different numbers of modes, $N$. Proposition~\ref{prop:HinfDesign} guarantees $J_5(t)\le 0$ (blue line). However, if we include one more mode without adjusting $\gamma$ and the controller, {then 
\begin{equation*}
    J_6(t)=J_5(t)+\int_0^t\left[(1+\rho_x\omega_6^2)z_6^2(t)-\gamma^2w_6^2(t)\right]\,dt
\end{equation*}
becomes positive for $t>20$ (green line). This happens because the $L^2$ gain for the additional mode with $n=6$ is greater than $\gamma$ and the additional integral term is positive.} The red line shows $J_{50}(t)\approx J_\infty(t)${, which is the cost when all the modes are considered.} Clearly, the controller designed using only $5$ modes cannot guarantee \eqref{HinfObjective} for the original system. 

Spillover occurs because the effect of the controller on the truncated modes is ignored. In the remainder of the paper, we provide a simple remedy to avoid spillover. Namely, we show how to modify $\rho_u$ in \eqref{J(t)} so that a controller guaranteeing $J_N(t)\le0$ for a given $N$ will guarantee \eqref{HinfObjective} with the original $\rho_u$. 
\subsection{Cost decomposition}\label{sec:cost decomposition}
We represent \eqref{FourierCoefficients} as 
\begin{subequations}\label{ODE:MD}
\begin{align}
\dot z^N&=Az^N+Bu+Ew^N,\label{ODE:MDa} \\
\dot{\bar z}_n&=A_n\bar z_n+B_nu+E_nw_n,\quad n>N, \label{ODE:MDb}
\end{align}
\end{subequations}
where $N\in\N$,
\begin{equation}\label{FD definitions}
\begin{aligned}
    &{z^N=\left[\begin{smallmatrix}
    z_1 \\
    \vdots \\
    z_N \\
    \dot{z}_1 \\
    \vdots \\
    \dot{z}_N
    \end{smallmatrix}\right],}\ 
    w^N=\left[\begin{smallmatrix}
    w_1 \\
    \vdots \\
    w_N
    \end{smallmatrix}\right],\ 
    B=\left[\begin{smallmatrix}
        0 \\
        \vdots \\
        0 \\
        b_1 \\
        \vdots \\
        b_N
        \end{smallmatrix}\right],\ 
    E=\left[\begin{smallmatrix}
    0_N\\ I_N
    \end{smallmatrix}\right],\\
    &A=\left[\begin{smallmatrix}
    0_N & I_N \\
    -\Omega_N^2 & -\left(c_1 I_N+c_2 \Omega_N^2\right)
    \end{smallmatrix}\right], \quad \Omega_N=\operatorname{diag}\{\omega_1,\ldots,\omega_N\}, 
\end{aligned}
\end{equation}
and the remaining notations are from \eqref{FourierCoefficients}. We will design an $H_\infty$ controller for \eqref{ODE:MDa} with the cost, $J_0$, that accounts for its effect on \eqref{ODE:MDb}. To find this cost, we decompose the original cost $J$ from \eqref{HinfObjective}. Namely, since $z(\cdot,t)\in H^2(0,\pi)$ (see Section~\ref{sec:Well-posedness}), Parseval's identity gives
\begin{equation*}
    \|z(\cdot,t)\|^2=\sum_{n=1}^\infty z_n^2(t),\qquad 
    \|z_{xx}(\cdot,t)\|^2=\sum_{n=1}^\infty \omega_n^2z_n^2(t). 
\end{equation*}
Our key idea is to represent $J$ from \eqref{HinfObjective} as 
\begin{equation}\label{Jdecomp}
\textstyle J=J_0+\sum_{n=N+1}^\infty J_n, 
\end{equation}
where 
\begin{equation*}
{\begin{aligned}
J_0&\textstyle =\int_0^\infty\Bigl[\sum_{n=1}^N(1+\rho_x\omega_n^2)z_n^2(t)+\left(\rho_u+\sum_{n=N+1}^\infty\rho_n\right)u^2(t)\\
&\textstyle \hspace{5cm}-\gamma^2\sum_{n=1}^Nw_n^2(t)\Bigr]\,dt,\\
J_n&\textstyle=\int_0^\infty\left[(1+\rho_x\omega_n^2)z_n^2(t)-\rho_nu^2(t)-\gamma^2w_n^2(t)\right]\,dt.
\end{aligned}}
\end{equation*}
The control, $u(t)$, is treated as a disturbance in \eqref{ODE:MDb}. Using the bounded real lemma (Corollary~\ref{cor:BRL}), we will find the minimum $\rho_n$ such that $J_n\le0$ for the zero initial conditions and any $w_n\in L^2([0,\infty),\R)$. Then, we will show that  $\sum_{n=N+1}^\infty\rho_n<\infty$ and construct a controller for \eqref{ODE:MDa} guaranteeing $J_0\le0$. 
\subsection{Bounded real lemma for the residue}
For a given $n>N$, \eqref{ODE:MDb} can be represented as \eqref{LTI} with 
\begin{equation*}
\begin{aligned}
&x=\bar z_n,\quad A=A_n,\quad B=0_{2\times 1},\\
&v=\left[\begin{smallmatrix}
\sqrt{\rho_n}u/\gamma \\ w_n
\end{smallmatrix}\right], \quad E=\begin{bmatrix}\frac{\gamma}{\sqrt{\rho_n}}B_n & E_n\end{bmatrix}. 
\end{aligned}
\end{equation*}
Note that the control input, $u$, is considered as a part of the disturbance, $v$, since the $H_\infty$ control will be designed based on \eqref{ODE:MDa}. The cost in \eqref{HinfProblem} coincides with $J_n$ for 
\begin{equation*}
C=\begin{bmatrix}
\sqrt{1+\rho_x\omega_n^2} & 0 
\end{bmatrix}\quad\text{and}\quad D=0. 
\end{equation*}
Then, the algebraic Riccati equation \eqref{BRL:ARE} takes the form 
\begin{equation}\label{ARE for n}
P_nA_n+A_n^\top P_n+\gamma^{-2}P_n\left[\begin{smallmatrix}
0 & 0 \\ 0 & 1+\gamma^2b_n^2/\rho_n
\end{smallmatrix}\right]P_n+\left[\begin{smallmatrix}
1+\rho_x\omega_n^2 & 0 \\ 0 & 0
\end{smallmatrix}\right]=0. 
\end{equation}
In \ref{appendix:ARE solution}, we show that the smallest $\rho_n$ guaranteeing the feasibility of \eqref{ARE for n} is 
\begin{equation*}
    \rho_n=\left\{
    \begin{aligned}
    &\frac{b_n^2(1+\rho_x\omega_n^2)}{4\omega_n^4\zeta_n^2(1-\zeta_n^2)-(1+\rho_x\omega_n^2)\gamma^{-2}}&&\text{if }2\zeta_n^2\le1,\\
    &\frac{b_n^2(1+\rho_x\omega_n^2)}{\omega_n^4-(1+\rho_x\omega_n^2)\gamma^{-2}}&&\text{if }2\zeta_n^2>1. 
    \end{aligned}
    \right.
\end{equation*}
The value of $\rho_n$ is the $L^2$ gain from $u$ to $\bar z_n$. Corollary~\ref{cor:BRL} guarantees $J_n\le0$ for these $\rho_n$. This can be used to obtain the $L^2$ gain of \eqref{EB_beam} without control. 
\begin{prop}[$L^2$ gain without control]\label{prop:gamma0}
The $L^2$ gain of the control-free \eqref{EB_beam} subject to \eqref{c1c2cond} is not greater than 
\begin{equation*}
    \gamma_0=\frac{2\sqrt{1+\rho_x}}{(c_1+c_2)\sqrt{4-(c_1+c_2)^2}}. 
\end{equation*}
\end{prop}
{\em Proof.} Repeating the arguments of \ref{appendix:ARE solution} with $\alpha_n=\gamma^{-2}$, we obtain that \eqref{ARE for n} is feasible for any $n\in\N$ if 
\begin{subequations}\label{gammaNoControl}
\begin{align}
&\gamma^2\ge\frac{1+\rho_x\omega_n^2}{4\omega_n^4\zeta_n^2(1-\zeta_n^2)}&\text{when}\quad 2\zeta_n^2\le1,\label{gammaNoControl1}\\
&\gamma^2\ge\frac{1+\rho_x\omega_n^2}{\omega_n^4}&\text{when}\quad 2\zeta_n^2>1. \label{gammaNoControl2}
\end{align}
\end{subequations}
The right-hand sides of \eqref{gammaNoControl} are decreasing in $n$. Moreover, 
\begin{equation*}
    \omega_n^4-4\omega_n^4\zeta_n^2(1-\zeta_n^2)=\omega_n^4(1-2\zeta_n^2)^2\ge0
\end{equation*}
implies that the bound in \eqref{gammaNoControl1} is not smaller than in \eqref{gammaNoControl2}. Since \eqref{c1c2cond} guarantees $2\zeta_1^2\le1$, the lower bound on $\gamma$ is obtained from \eqref{gammaNoControl1} with $n=1$, i.e., with $\omega_1=1$ and $\zeta_1=(c_1+c_2)/2$. The feasibility of \eqref{ARE for n} implies $J_n\le0$. Taking $N=0$ and $J_0=0$ in \eqref{Jdecomp}, we obtain $J\le0$. 
\hfill$\square$

In \ref{appendix:rhon}, we show that 
\begin{equation}\label{rhoInfBound}
\sum_{n=N+1}^\infty\rho_n\le\rho_\infty=\sum_{n=N+1}^M\rho_n+C_M\left[|x_R-x_L|-\sum_{n=1}^M\frac{b_n^2}{\omega_n^2}\right], 
\end{equation}
where 
\begin{equation*}
\begin{aligned}
&C_M=\frac{\omega_{M+1}^2(1+\rho_x\omega_{M+1}^2)}{\omega_{M+1}^4-(1+\rho_x\omega_{M+1}^2)\gamma^{-2}},\\
&M=\max\left\{N,\left\lfloor\sqrt{\frac{1+\sqrt{1-2 c_1 c_2}}{\sqrt{2}c_2}}\right\rfloor\right\}. 
\end{aligned}
\end{equation*}
Here, $\lfloor\cdot\rfloor$ stands for the integer part. Note that \eqref{c1c2cond} implies $2c_1c_2\le1$. 

As explained in \ref{appendix:rhon}, $\sum_{n=1}^\infty b_n^2/\omega_n^2=|x_R-x_L|$. Therefore, $\rho_\infty\to 0$ monotonically as $N\to\infty$. That is, by considering more modes in the control design, we reduce the $L^2$ gain of the residue associated with the spillover. 
\subsection{$H_\infty$ controller design without spillover}\label{sec:HinfDesign}
The system \eqref{ODE:MDa} is in the form of \eqref{LTI} with $x=z^N$, $v=w^N$, and $A$, $B$, and $E$ defined in \eqref{FD definitions}. Taking 
\begin{equation}\label{CD}
    C=\left[\begin{smallmatrix}
    \sqrt{I_N+\rho_x\Omega_N^2} & 0_{N\times N} \\
    0_{1\times N} & 0_{1\times N}
    \end{smallmatrix}\right]\quad\text{and}\quad
    D=\left[\begin{smallmatrix}
    0_{N\times 1}\\\sqrt{\rho_u+\rho_\infty}
    \end{smallmatrix}\right], 
\end{equation}
we obtain that $D^\top C=0$, $R=D^\top D=\rho_u+\rho_\infty>0$, and the left-hand side of \eqref{HinfProblem} coincides with $J_0$ from \eqref{Jdecomp}. By Proposition~\ref{prop:HinfDesign}, if $0<P\in\R^{2N\times 2N}$ satisfies \eqref{ARE}, then 
\begin{equation}\label{feedback}
    u(t)=-(\rho_u+\rho_\infty)^{-1}B^\top Pz^N(t) 
\end{equation}
guarantees $J_0\le0$. Since $\rho_n$ were chosen so that $J_n\le0$, we obtain that $J\le0$. 

{Note that \eqref{w condition} implies $w^N\in L^2(0,\infty)$. Therefore, the solution of the stable system \eqref{ODE:MDa}, \eqref{feedback} satisfies $z^N\in L^2(0,\infty)$. That is, the right-hand side of \eqref{ODE:MDa} is from $L^2$, meaning that $\dot z^N\in L^2(0,\infty)$, and both $z^N\in H^1(0,\infty)$ and $u\in H^1(0,\infty)$. This is the property we used in the well-posedness analysis of Section~\ref{sec:Well-posedness}. Summarizing, we have the following result.} 

\begin{thm}\label{th:L2gain}
Consider the Euler--Bernoulli beam \eqref{EB_beam} subject to \eqref{c1c2cond} and its modal decomposition \eqref{ODE:MD} with some $N\in\N$. Given non-negative $\rho_x$, $\rho_u$, and $\gamma$, let $\rho_\infty$ be given by \eqref{rhoInfBound}. If $0<P\in\R^{2N\times 2N}$ satisfies the algebraic Riccati equation \eqref{ARE} with $A$, $B$, $C$, $D$, and $E$ given in \eqref{FD definitions} and \eqref{CD}, then the state feedback \eqref{feedback} guarantees that the $L^2$ gain of \eqref{EB_beam} is not greater than $\gamma$, that is, \eqref{HinfObjective} holds for $z(\cdot,0)\equiv0\equiv z_t(\cdot,0)$ and any $w$ satisfying \eqref{w condition}. 
\end{thm}
{
{\em Proof.} Let $M$ be as defined below (22). Consider 
\begin{equation*}
\textstyle V=(z^N)^\top Pz^N+\sum_{n=N+1}^\infty\bar z_n^\top P_n\bar z_n 
\end{equation*}
with $P_n$ defined by \eqref{Pn2} for $N+1\le n\le M$ and by \eqref{Pn1} for $n>M$. The series converges since 
\begin{equation}\label{Pnlimit}
    P_n\sim\rho_x\left[\begin{smallmatrix}
    \omega_n^2c_2 & 1 \\ 1 & 2c_2
    \end{smallmatrix}\right]\quad\text{as}\quad n\to\infty
\end{equation}
(we chose ``$+$'' for the right bottom element), while $z(\cdot,t)\in H^2(0,\pi)$ and $z_t(\cdot,t)\in L^2(0,\pi)$ (see Section~\ref{sec:Well-posedness}). Let $J(t)$ be $J$ as defined in \eqref{HinfObjective} but with $\infty$ replaced by $t$. Calculating the derivative along the trajectories of \eqref{ODE:MD} and using the relation $\rho_u+\sum_{N+1}^\infty\rho_n\le\rho_u+\rho_\infty=R$, we obtain 
\begin{equation*}
\begin{aligned}
\dot V(t)+\dot J(t)\le{}&2(z^N)^\top P[Az^N+Bu+Ew^N]\\
&+(z^N)^\top C^\top Cz^N+Ru^2-\gamma^2|w^N|^2\\
&\textstyle{}+2\sum_{n=N+1}^\infty \bar z_n^\top P_n[A_n\bar z_n+B_nu+E_nw_n]\\
&\textstyle{}+\sum_{n=N+1}^\infty[(1+\rho_x\omega_n^2)z_n^2-\rho_nu^2-\gamma^2w_n^2]. 
\end{aligned}
\end{equation*}
Completing the squares, we find 
\begin{equation}\label{complitionOfSquares}
\begin{aligned}
&2(z^N)^\top PBu+Ru^2=|R^{-\frac12}B^\top Pz^N+R^{\frac12}u|^2-(z^N)^\top PBR^{-1}B^\top Pz^N, \\
&2(z^N)^\top PEw^N-\gamma^2|w^N|^2=\gamma^{-2}|E^\top Pz^N|^2-|\gamma^{-1}E^\top Pz^N-\gamma w^N|^2,\\
&2\bar z_n^\top P_nB_nu-\rho_nu^2=\rho_n^{-1}|B_n^\top P_n\bar z_n|^2-|\rho_n^{-\frac12}B_n^\top P_n\bar z_n-\rho_n^{\frac12}u|^2,\\
&2\bar z_n^\top P_nE_nw_n-\gamma^2w_n^2=\gamma^{-2}|E_n^\top P_n\bar z_n|^2-|\gamma^{-1}E_n^\top P_n\bar z_n-\gamma w_n|^2. 
\end{aligned}
\end{equation}
In view of \eqref{ARE} and \eqref{ARE for n}, these lead to 
\begin{multline*}
\dot V(t)+\dot J(t)\le|R^{-\frac12}B^\top Pz^N+R^{\frac12}u|^2-|\gamma^{-1}E^\top Pz^N-\gamma w^N|^2\\
-\sum_{n=N+1}^\infty\left[|\rho_n^{-\frac12}B_n^\top P_n\bar z_n-\rho_n^{\frac12}u|^2+|\gamma^{-1}E_n^\top P_n\bar z_n-\gamma w_n|^2\right]. 
\end{multline*}
Substituting $u$ from \eqref{feedback}, we obtain 
\begin{equation}\label{dotVdotJ}
    \dot V(t)+\dot J(t)\le 0. 
\end{equation}
Integrating the above from $0$ to $t$, we obtain 
\begin{equation*}
    V(t)-V(0)+J(t)-J(0)\le 0. 
\end{equation*}
Given that $V(0)=0$ for the zero initial conditions, and $J(0)=0$, we have $J(t)\le -V(t)\le 0$, which implies \eqref{HinfObjective}. \qed}

\begin{rem}[Internal stability]
The designed feedback \eqref{feedback} renders \eqref{EB_beam} internally stable in the norm 
\begin{equation*}
    \|z(\cdot,t)\|_X^2=\|z_{xx}(\cdot,t)\|^2+\|z_t(\cdot,t)\|^2. 
\end{equation*}
Indeed, \eqref{Pnlimit} implies the existence of positive $\varepsilon_1$ and $\varepsilon_2$ such that $\varepsilon_1\|z(\cdot,t)\|_X^2\le V\le \varepsilon_2\|z(\cdot,t)\|_X^2$, and \eqref{dotVdotJ} implies $\dot V\le0$ for $w(\cdot,t)\equiv0$. 
\end{rem}

\begin{rem}[Solution existence]
Since $A$, defined in \eqref{FD definitions}, is Hurwitz, $(A,B)$ is stabilizable. It is easy to check that $(A,C)$ is observable, hence detectable. As mentioned in Remark~\ref{rem:AREexistence}, this guarantees that \eqref{ARE} has a solution for a large enough $\gamma$. That is, the conditions of Theorem~\ref{th:L2gain} hold for any $N\in\N$ and large enough $\gamma$. 
\end{rem}
\begin{rem}[Number of modes and the $L^2$ gain]\label{rem:N vs gamma}
When $N$ grows, $\gamma$ can only decrease. Indeed, we know that $J_{N+1}(t)\le0$ and \eqref{feedback} guarantees $J_0(t)\le0$ with $J_{N+1}(t)$ and $J_0(t)$ defined below \eqref{Jdecomp}. Taking $K_1,K_2\in\R^{1\times N}$ such that $[K_1\ K_2]=(\rho_u+\rho_\infty)^{-1}B^\top P$, we have that 
\begin{equation*}
    u=-\begin{bmatrix}K_1 & 0 & K_2 & 0 \end{bmatrix}z^{N+1}
\end{equation*} 
guarantees $\bar J_0(t)=J_0(t)+J_{N+1}(t)\le0$. Note that $\bar J_0(t)$ is $J_0(t)$ with $N$ replaced by $N+1$. By \cite[Theorem~6.3.6]{Green2012}, \eqref{ARE} has a solution for the matrices defined in \eqref{FD definitions} and \eqref{CD} with $N$ replaced by $N+1$. That is, the same $\gamma$ is achievable with $N+1$ modes. When considering $N+1$ modes, we are making the sum $J_0(t)+J_{N+1}(t)$ negative instead of each term, $J_0(t)$ and $J_{N+1}(t)$, independently. This gives more flexibility and may reduce $\gamma$, as demonstrated in Fig.~\ref{fig:gamma_vs_N}. 
\end{rem}



\section{Numerical simulations}\label{sec:example}
\begin{figure}
	\centering
	\includegraphics[width=\linewidth]{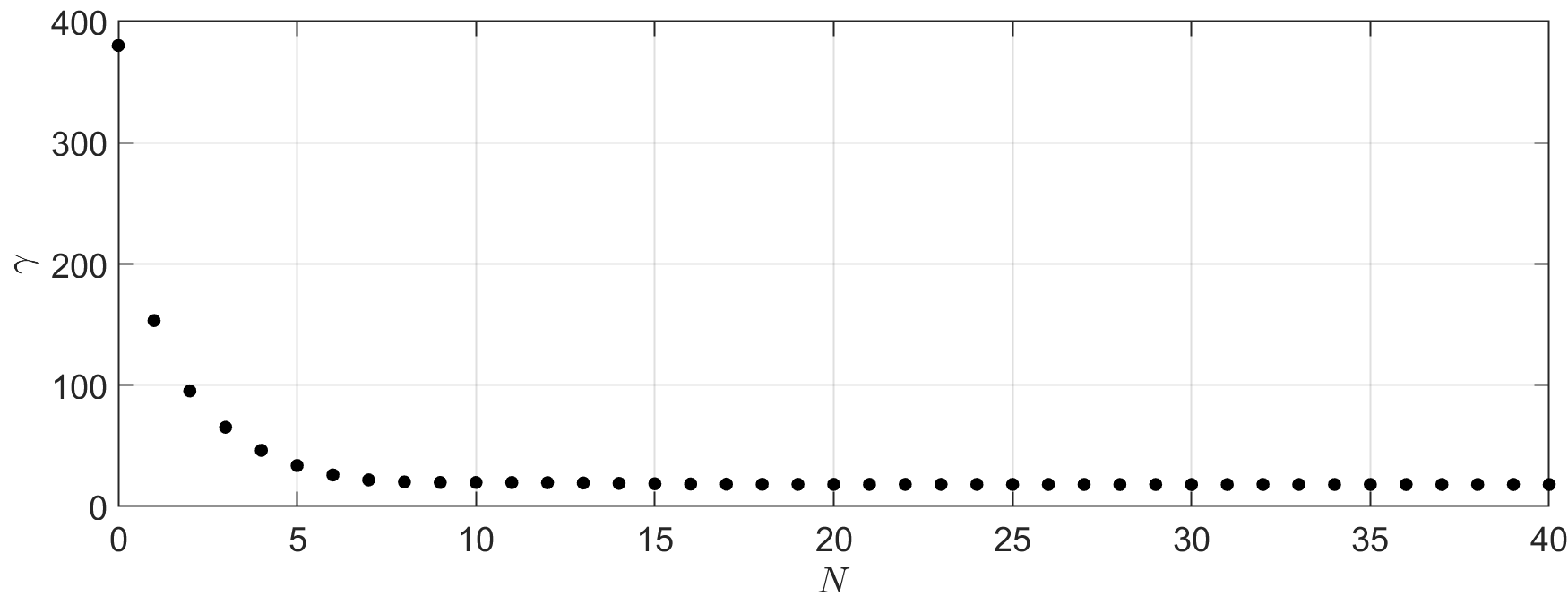}
	\caption{The $L^2$ gain of the Euler--Bernoulli beam \eqref{EB_beam} for different numbers of controlled modes $N$.}\label{fig:gamma_vs_N}
\end{figure}

As an example, we consider an aluminum rectangular beam of dimensions $1\:m\times0.1\:m\times0.01\:m$ with hinged ends and a piezoelectric actuator of length $2\:cm$ placed at $30\:cm$ from the left edge. This system can be modeled by \eqref{EB_beam_original} with the parameters given in the following table: 
\begin{center}
\begin{tabular}{lll}
\hline
Linear density & $\mu$ & $2.71\:kg/m$\\
Young's modulus & $E$ & $70\times 10^9\:N/m^2$\\
Moment of inertia & $I$ & $8.3\times 10^{-8}\:m^4$\\
Viscous damping & $c_v$ & $1.76\:kg/(m\cdot s)$\\
Structural damping & $c_k$ & $2.05\times 10^5\:kg/(m\cdot s)$\\
Left actuator position & $\tilde x_L$ & $0.29\:m$\\
Right actuator position & $\tilde x_R$ & $0.31\:m$\\
\hline
\end{tabular}
\end{center}
The linear density is calculated as $\mu=\rho A$, where $\rho=2710\: kg/m^3$ is the density of aluminum, and $A=0.1\times 0.01=10^{-3}\:m^2$ is the cross-section area of the beam. The damping coefficients, $c_v$ and $c_k$, are taken from \cite{Banks1991}. The value of $c_a$ depends on the type of the piezoelectric patch; it does not affect the performance analysis since the control can be scaled as $\tilde u'=c_a\tilde u$. After the change of variables \eqref{changeOfVariables}, we obtain \eqref{EB_beam} with 
\begin{equation*}
    c_1=1.4\times 10^{-3},\ c_2=1.3\times 10^{-3},\ x_L=0.91,\ x_R=0.97. 
\end{equation*}
Our objective is to design a state-feedback control law of the form \eqref{feedback} guaranteeing that the solution of \eqref{EB_beam} with $z(\cdot,0)\equiv0\equiv z_t(\cdot,0)$ satisfies \eqref{HinfObjective} with $\rho_u=10^{-3}$, $\rho_x=0.1$, and smallest possible $\gamma>0$. To decide on how many modes to consider in the controller design, we calculate the minimum $\gamma$ for different numbers of controlled modes,~$N$. Proposition~\ref{prop:gamma0} gives $\gamma_0\approx 380$ as the smallest $L^2$ gain without control. For each integer $N\in[1,40]$, we found the minimum $\gamma$ satisfying the conditions of Theorem~\ref{th:L2gain}. The results are shown in Fig.~\ref{fig:gamma_vs_N}. As explained in Remark~\ref{rem:N vs gamma}, the $L^2$ gain decreases when more modes are considered. The limit value is $\gamma\approx 18$. Since $\gamma$ does not improve significantly for $N>8$, we consider $N=8$ modes. In this case, $\gamma\approx20.2$ and $\rho_\infty\approx8\times 10^{-3}$, which we found using \eqref{rhoInfBound}. To find the controller gain in \eqref{feedback}, we solve \eqref{ARE} for $P>0$ with $A$, $B$, $C$, $D$, and $E$ defined in \eqref{FD definitions} and \eqref{CD}. Note that, for this example, the first condition in (15) of \cite{Selivanov2023} requires $N\ge32$ and the resulting LMIs are not feasible for any $\gamma>0$. 

\begin{figure}
	\centering
	\includegraphics[width=\linewidth]{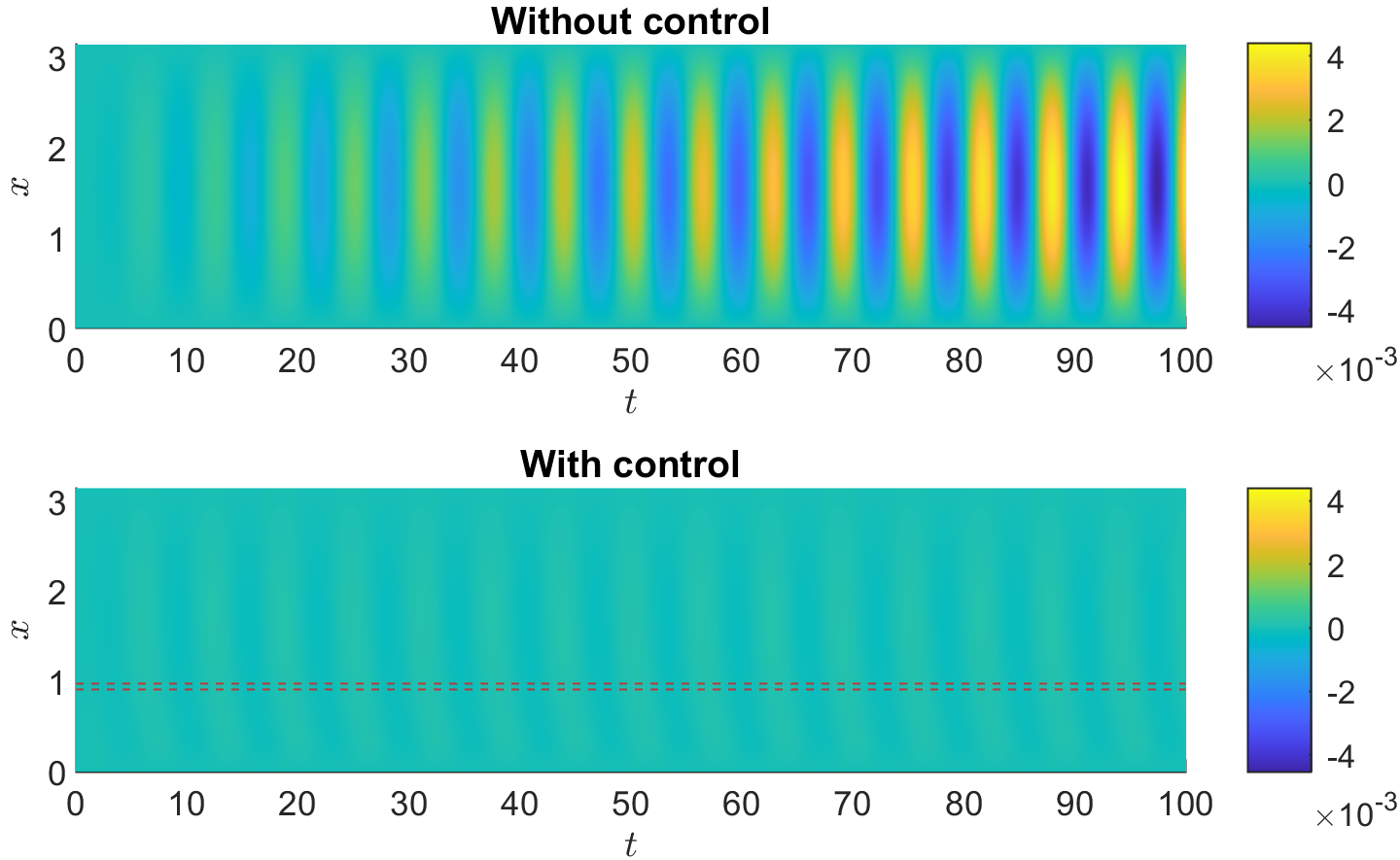}
	\caption{Euler--Bernoulli beam without and with control. The red dashed lines show the ends of the piezoelectric actuator.}\label{fig:States}
\end{figure}

\begin{figure}
	\centering
	\includegraphics[width=\linewidth]{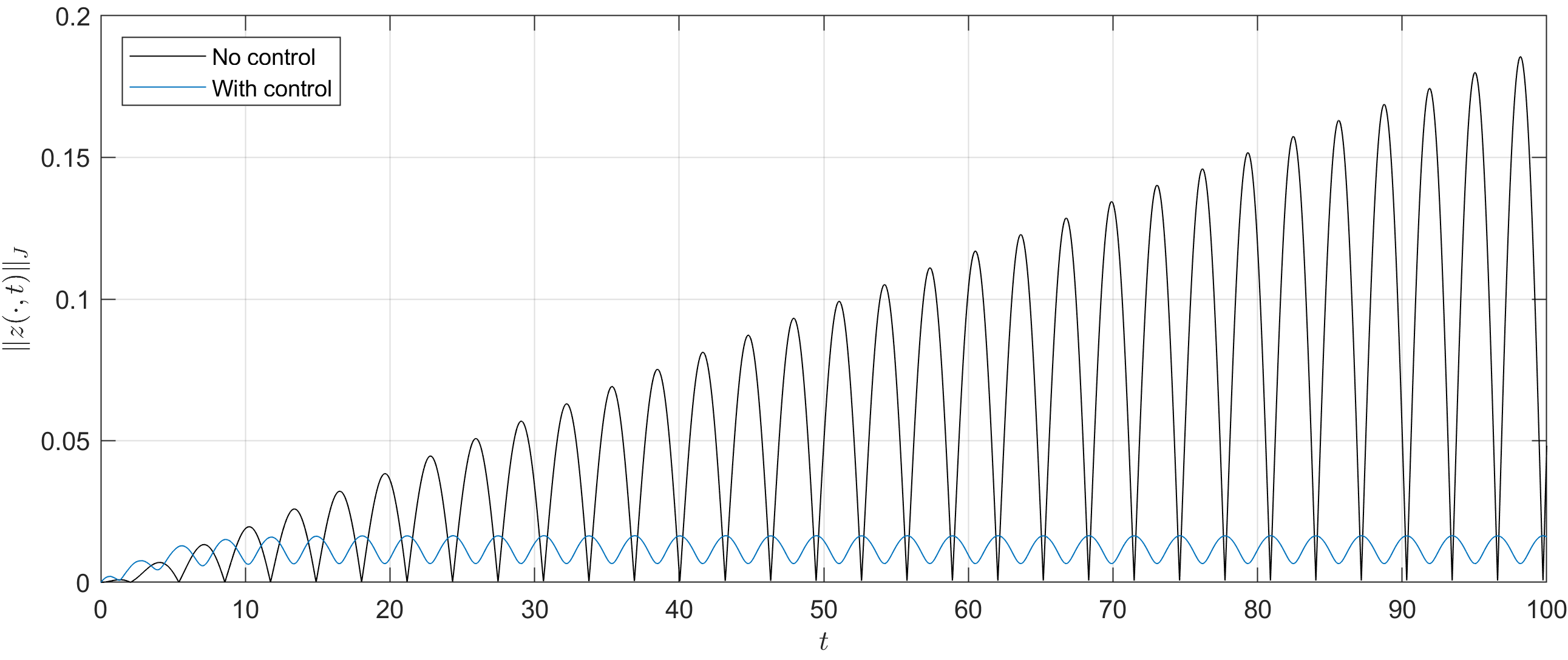}
	\caption{The value of $\|z(\cdot,t)\|_J$, defined in \eqref{Jnorm}, without (black) and with (blue) control.}\label{fig:H2norm}
\end{figure}

The results of numerical simulations without and with control for the same disturbance are shown in Fig.~\ref{fig:States}. {To generate the disturbance,} we calculated $P>0$ satisfying \eqref{BRL:ARE} with $A$, $E$, and $C$ given in \eqref{FD definitions} and \eqref{CD}, found $z_d^N(t)$ as the solution of \eqref{ODE:MDa} with $N=30$, $u\equiv 0$, and $z_d^N(0)=[1,\ldots,1]^\top\in\R^{60}$, substituted $w^N(t)=\gamma^{-2}E^\top Pz_d^N(t)$ into \eqref{FD definitions}, and took $w(x,t)=\sum_{n=1}^Nw_n(t)\varphi_n(x)$. {The value of $w^N$ was selected to maximizes the related negative term in \eqref{complitionOfSquares}.} Clearly, the proposed control strategy attenuates the effect of the disturbance. This is also evident from Fig.~\ref{fig:H2norm}, which shows 
\begin{equation}\label{Jnorm}
    \|z(\cdot,t)\|_J=\sqrt{\|z(\cdot,t)\|^2+\rho_x\|z_{xx}(\cdot,t)\|^2}
\end{equation}
without (black) and with (blue) control. 
\begin{figure}
	\centering
	\includegraphics[width=\linewidth]{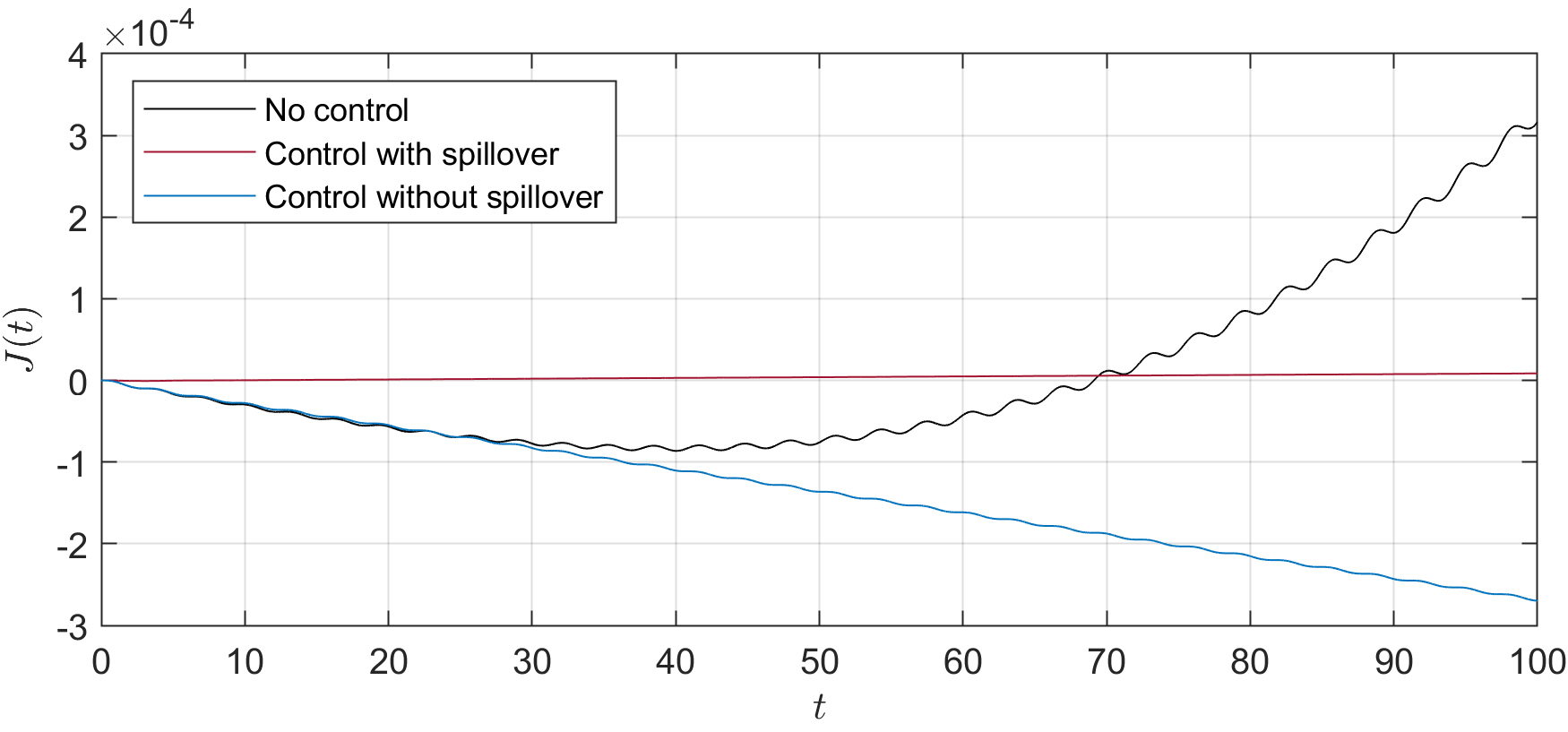}
	\caption{The value of $J(t)$ (given by \eqref{HinfObjective} with $\infty$ replaced by $t$) without control (black) and with control (blue) for $\gamma\approx 20.2$. The red line shows the spillover phenomenon occurring when the modes with $n>8$ are ignored.}\label{fig:J}
\end{figure}

The value of $J(t)$, obtained by replacing $\infty$ with $t$ in \eqref{HinfObjective}, is shown in Fig.~\ref{fig:J}. As guaranteed by Theorem~\ref{th:L2gain}, the control ensures that $J=\lim_{t\to\infty}J(t)<0$ for $\gamma\approx 20.2$ (blue line). Without control (black line), $J(t)$ becomes positive for $t\approx 70$. If the residue is ignored ($\rho_\infty=0$), a smaller $\gamma\approx 7.16$ is obtained following the steps detailed in Section~\ref{sec:spillover}. In this case, the spillover phenomenon causes $J(t)>0$ for $t>7$ (red line). This vividly demonstrates why the residue, i.e., the modes with $n>N$, must not be ignored. Theorem~\ref{th:L2gain} provides a simple way of designing a controller avoiding the spillover phenomenon. 
\begin{rem}[$N$ vs $\gamma$ with spillover]
If the residue is ignored ($\rho_\infty=0$), then $\gamma\approx 6.97$ for $N=5$ (see Section~\ref{sec:spillover}) and $\gamma\approx 7.16$ for $N=8$ (as explained above). That is, the $L^2$ gain may increase when more modes are considered. This happens because, by increasing $N$, one obtains a more accurate estimate of the actual $L^2$ gain, which is larger than that obtained using the truncated modal decomposition. If the residue is accounted for, larger $N$ will never lead to a larger $\gamma$ (see Remark~\ref{rem:N vs gamma}). 
\end{rem}
\section{Conclusions}
We studied the $H_\infty$ control of the Euler--Bernoulli beam with viscous and Kelvin--Voigt damping using piezoelectric actuators. We showed that spillover occurs when a finite number of modes are considered in the $H_\infty$ design. Then we proposed a simple modification of the cost guaranteeing that the controller designed based on a finite number of modes does not lead to spillover. Using a realistic model of the beam, we demonstrated how to find the number of modes required to design a controller, i.e., such that a further increase of the number of considered modes does not improve the $L^2$ gain significantly. 

\section*{Acknowledgements}
The work of Emilia Fridman was supported in part by the ISF-NSFC Joint Research Program under Grant 3054/23 and in part by the C. and H. Manderman Chair at Tel Aviv University.

\providecommand{\noopsort}[1]{}

\appendix
\section{Solution of \eqref{ARE for n}}\label{appendix:ARE solution}
Let $P_n=\left[\begin{smallmatrix} p_1 & p_2 \\ p_2 & p_3 \end{smallmatrix}\right]$. Then \eqref{ARE for n} is equivalent to 
\begin{equation*}
\begin{aligned}
&\alpha_np_2^2-2\omega_n^2p_2+(1+\rho_x\omega_n^2)=0,\\
&p_1-2\zeta_n\omega_np_2-p_3\omega_n^2+\alpha_np_2p_3=0,\\
&\alpha_np_3^2-4\zeta_n\omega_np_3+2p_2=0 
\end{aligned}
\end{equation*}
with $\alpha_n=b_n^2/\rho_n+\gamma^{-2}$. These are equivalent to 
\begin{equation}\label{p1:3}
\begin{aligned}
&p_2=\alpha_n^{-1}\left[\omega_n^2\pm\sqrt{\omega_n^4-\alpha_n(1+\rho_x\omega_n^2)}\right],\\
&p_3=\alpha_n^{-1}\left[2\zeta_n\omega_n \pm \sqrt{4\zeta_n^2\omega_n^2-2\alpha_np_2}\right],\\
&p_1=2\zeta_n\omega_np_2+p_3\omega_n^2-\alpha_np_2p_3.
\end{aligned}
\end{equation}
These values are real if and only if 
\begin{subequations}
\begin{align}
&\omega_n^4\ge\alpha_n(1+\rho_x\omega_n^2)\quad\text{and}\label{alphan:p2cond}\\
&2\zeta_n^2\omega_n^2\ge\alpha_np_2=\omega_n^2\!-\!\sqrt{\omega_n^4\!-\!\alpha_n(1\!+\!\rho_x\omega_n^2)}. \label{alphan:p3cond}
\end{align}
\end{subequations}
We took $p_2$ with ``$-$'' since $2\zeta_n^2$ can be smaller than $1$. 

To minimize $\rho_n$, we maximize $\alpha_n$. {If $2\zeta_n^2<1$, then \eqref{alphan:p3cond} gives the maximum $\alpha_n=4\omega_n^4\zeta_n^2(1-\zeta_n^2)/(1+\rho_x\omega_n^2)$, which satisfies \eqref{alphan:p2cond} since $\omega_n^4-\alpha_n(1+\rho_x\omega_n^2)=\omega_n^4(1-2\zeta_n^2)^2>0$. Substituting this into \eqref{p1:3}, we obtain 
\begin{equation}\label{Pn2}
    P_n=\frac{2\zeta_n\omega_n}{\alpha_n}\left[\begin{smallmatrix}
    \omega_n^2 & 
    \zeta_n\omega_n\\
    \zeta_n\omega_n & 1
    \end{smallmatrix}\right]>0. 
\end{equation}
If $2\zeta_n^2\ge1$, then \eqref{alphan:p3cond} is true subject to \eqref{alphan:p2cond}, which gives $\alpha_n=\omega_n^4/(1+\rho_x\omega_n^2)$. Substituting this into \eqref{p1:3}, we obtain 
\begin{equation}\label{Pn1}
    P_n=\frac{\omega_n}{\alpha_n}
    \left[\begin{smallmatrix}
    2\zeta_n\omega_n^2 & 
    \omega_n \\
    \omega_n  & 2\zeta_n\pm\sqrt{4\zeta_n^2-2}
    \end{smallmatrix}\right]>0. 
\end{equation}}%
The minimum values of $\rho_n$ are calculated from $\alpha_n=b_n^2/\rho_n+\gamma^{-2}$ with the corresponding $\alpha_n$. 
\section{Upper bound on the $L^2$ gain for the residue}\label{appendix:rhon}
For $n>M$, we have $\zeta_n\ge1/\sqrt{2}$. Then 
\begin{equation*}
\begin{aligned}
\rho_n&\textstyle=\frac{b_n^2(1+\rho_x\omega_n^2)}{\omega_n^4-(1+\rho_x\omega_n^2)\gamma^{-2}}=\frac{\omega_n^{-2}+\rho_x}{1-(\omega_n^{-2}+\rho_x)\omega_n^{-2}\gamma^{-2}}\frac{b_n^2}{\omega_n^2}\\
&\textstyle\le\frac{\omega_{M+1}^{-2}+\rho_x}{1-(\omega_{M+1}^{-2}+\rho_x)\omega_{M+1}^{-2}\gamma^{-2}}\frac{b_n^2}{\omega_n^2}=C_M\frac{b_n^2}{\omega_n^2}. 
\end{aligned}
\end{equation*}
Note that $b_n/\omega_n$ are the Fourier coefficients of 
\begin{equation*}
    \chi_{[x_L,x_R]}(x)=\left\{
    \begin{aligned}
    &1,&& x\in[x_L,x_R],\\
    &0,&& x\notin[x_L,x_R].
    \end{aligned}
    \right.
\end{equation*}
By Parseval's identity, 
\begin{equation*}
\textstyle\sum_{n=1}^\infty\frac{b_n^2}{\omega_n^2}=\|\chi_{[x_L,x_R]}\|^2=|x_R-x_L|. 
\end{equation*}
Therefore, 
\begin{equation*}
    \sum_{n=M+1}^\infty\rho_n\le C_M\sum_{n=M+1}^\infty\frac{b_n^2}{\omega_n^2}=C_M\left[|x_R-x_L|-\sum_{n=1}^M\frac{b_n^2}{\omega_n^2}\right], 
\end{equation*}
which implies \eqref{rhoInfBound}. 

\begin{thebibliography}{10}

\bibitem{Francis1987}
B.~A. Francis, {\em A {{Course}} in {{$H_\infty$ Control Theory}}},
  vol.~88 of {\em Lecture {{Notes}} in {{Control}} and {{Information
  Sciences}}}.
\newblock Springer, 1987.

\bibitem{Green2012}
M.~Green and D.~Limebeer, {\em Linear Robust Control}.
\newblock Dover, 2012.

\bibitem{ozbay1990}
H.~Özbay and A.~Tannenbaum, ``A {{Skew Toeplitz Approach}} to the {{$H^\infty$
  Optimal Control}} of {{Multivariable Distributed Systems}},'' {\em SIAM
  Journal on Control and Optimization}, vol.~28, no.~3, pp.~653--670, 1990.

\bibitem{ozbay1993a}
H.~Özbay, ``$H^\infty$ optimal controller design for a class of distributed
  parameter systems,'' {\em International Journal of Control}, vol.~58, no.~4,
  pp.~739--782, 1993.

\bibitem{Nett1983}
C.~N. Nett, C.~A. Jacobson, and M.~J. Balas, ``Fractional representation
  theory: {{Robustness}} results with applications to finite dimensional
  control of a class of linear distributed systems,'' in {\em The 22nd {{IEEE
  Conference}} on {{Decision}} and {{Control}}}, pp.~268--280, 1983.

\bibitem{Glover1986}
K.~Glover, ``Robust stabilization of linear multivariable systems: Relations to
  approximation,'' {\em International Journal of Control}, vol.~43, no.~3,
  pp.~741--766, 1986.

\bibitem{ozbay1991}
H.~Özbay and A.~Tannenbaum, ``On the structure of suboptimal {{$H^\infty$}}
  controllers in the sensitivity minimization problem for distributed stable
  plants,'' {\em Automatica}, vol.~27, no.~2, pp.~293--305, 1991.

\bibitem{Bontsema1988}
J.~Bontsema and R.~F. Curtain, ``A {{Note}} on {{Spillover}} and {{Robustness}}
  for {{Flexible Systems}},'' {\em IEEE Transactions on Automatic Control},
  vol.~33, no.~6, pp.~567--569, 1988.

\bibitem{Curtain1995}
R.~F. Curtain and H.~Zwart, {\em An {{Introduction}} to {{Infinite-Dimensional
  Linear Systems Theory}}}. Springer, 1995.

\bibitem{Curtain1990a}
R.~Curtain, ``$H_\infty$ control for distributed parameter systems: A survey,''
  in {\em 29th {{IEEE Conference}} on {{Decision}} and {{Control}}},
  pp.~22--26, 1990.

\bibitem{Keulen1993}
B.~{\noopsort{keulen}}van Keulen, {\em $H_\infty$-{{Control}} for {{Distributed
  Parameter Systems}}: {{A State-Space Approach}}}.
\newblock Birkhäuser, 1993.

\bibitem{Athans1970}
M.~Athans, ``Toward a practical theory for distributed parameter systems,''
  {\em IEEE Transactions on Automatic Control}, vol.~15, no.~2, pp.~245--247,
  1970.

\bibitem{Balas1982a}
M.~J. Balas, ``Toward {{A More Practical Control Theory}} for {{Distributed
  Parameter Systems}},'' in {\em Advances in {{Theory}} and {{Applications}}}
  (C.~Leondes, ed.), vol.~18, pp.~361--421, Academic Press, 1982.

\bibitem{Triggiani1980}
R.~Triggiani, ``Boundary {{Feedback Stabilizability}} of {{Parabolic
  Equations}},'' {\em Applied Mathematics and Optimization}, vol.~6, no.~1,
  pp.~201--220, 1980.

\bibitem{Balas1982}
M.~J. Balas, ``Trends in {{Large Space Structure Control Theory}}: {{Fondest
  Hopes}}, {{Wildest Dreams}},'' {\em IEEE Transactions on Automatic Control},
  vol.~27, no.~3, pp.~522--535, 1982.

\bibitem{Curtain1993}
R.~F. Curtain, ``A {{Comparison}} of finite-dimensional controller designs for
  distributed parameter systems,'' Research {{Report}} RR-1647, Inria, 1992.

\bibitem{Balas1978a}
M.~J. Balas, ``Feedback {{Control}} of {{Flexible Systems}},'' {\em IEEE
  Transactions on Automatic Control}, vol.~23, no.~4, pp.~673--679, 1978.

\bibitem{Meirovitch1983}
L.~Meirovitch and H.~Baruh, ``On the problem of observation spillover in
  self-adjoint distributed-parameter systems,'' {\em Journal of Optimization
  Theory and Applications}, vol.~39, no.~2, pp.~269--291, 1983.

\bibitem{Balas1988}
M.~J. Balas, ``Finite-dimensional controllers for linear distributed parameter
  systems: {{Exponential}} stability using residual mode filters,'' {\em
  Journal of Mathematical Analysis and Applications}, vol.~133, no.~2,
  pp.~283--296, 1988.

\bibitem{Moheimani1998}
S.~O.~R. Moheimani, ``Minimizing the {{Effect}} of {{Out}} of {{Bandwidth
  Modes}} in {{Truncated Structure Models}},'' {\em Journal of Dynamic Systems,
  Measurement, and Control}, vol.~122, no.~1, pp.~237--239, 1998.

\bibitem{Harkort2011}
C.~Harkort and J.~Deutscher, ``Finite-dimensional observer-based control of
  linear distributed parameter systems using cascaded output observers,'' {\em
  International Journal of Control}, vol.~84, no.~1, pp.~107--122, 2011.

\bibitem{Lasiecka1983}
I.~Lasiecka and R.~Triggiani, ``Stabilization and {{Structural Assignment}} of
  {{Dirichlet Boundary Feedback Parabolic Equations}},'' {\em SIAM Journal on
  Control and Optimization}, vol.~21, no.~5, pp.~766--803, 1983.

\bibitem{Curtain1984}
R.~F. Curtain, ``Finite {{Dimensional Compensators}} for {{Parabolic
  Distributed Systems}} with {{Unbounded Control}} and {{Observation}},'' {\em
  SIAM Journal on Control and Optimization}, vol.~22, no.~2, pp.~255--276,
  1984.

\bibitem{Karafyllis2016a}
I.~Karafyllis and M.~Krstic, ``{{ISS With Respect To Boundary Disturbances}}
  for 1-{{D Parabolic PDEs}},'' {\em IEEE Transactions on Automatic Control},
  vol.~61, no.~12, pp.~1--23, 2016.

\bibitem{Karafyllis2019k}
I.~Karafyllis and M.~Krstic, {\em Input-to-{{State Stability}} for {{PDEs}}}.
\newblock Springer, 2019.

\bibitem{Katz2020a}
R.~Katz and E.~Fridman, ``Constructive method for finite-dimensional
  observer-based control of 1-{{D}} parabolic {{PDEs}},'' {\em Automatica},
  vol.~122, p.~109285, 2020.

\bibitem{Karafyllis2021}
I.~Karafyllis, ``Lyapunov-based boundary feedback design for parabolic
  {{PDEs}},'' {\em International Journal of Control}, vol.~94, no.~5,
  pp.~1247--1260, 2021.

\bibitem{Selivanov2018e}
A.~Selivanov and E.~Fridman, ``Boundary {{Observers}} for a
  {{Reaction-Diffusion System Under Time-Delayed}} and {{Sampled-Data
  Measurements}},'' {\em IEEE Transactions on Automatic Control}, vol.~64,
  no.~8, pp.~3385--3390, 2019.

\bibitem{Lhachemi2021}
H.~Lhachemi and C.~Prieur, ``Predictor-based output feedback stabilization of
  an input delayed parabolic {{PDE}} with boundary measurement,'' {\em
  Automatica}, vol.~137, p.~110115, 2022.

\bibitem{Katz2022}
R.~Katz and E.~Fridman, ``Delayed finite-dimensional observer-based control of
  {{1D}} parabolic {{PDEs}} via reduced-order {{LMIs}},'' {\em Automatica},
  vol.~142, 2022.

\bibitem{Katz2021c}
R.~Katz and E.~Fridman, ``Global stabilization of a {{1D}} semilinear heat
  equation via modal decomposition and direct {{Lyapunov}} approach,'' {\em
  Automatica}, vol.~149, p.~110809, 2023.

\bibitem{Selivanov2024a}
A.~Selivanov, P.~Wang, and E.~Fridman, ``Guaranteed {{Cost Boundary Control}}
  of the {{Semilinear Heat Equation}},'' {\em IEEE Control Systems Letters},
  vol.~8, pp.~898--903, 2024.

\bibitem{Katz2022a}
R.~Katz and E.~Fridman, ``Finite-{{Dimensional Boundary Control}} of the
  {{Linear Kuramoto-Sivashinsky Equation Under Point Measurement With
  Guaranteed L2-Gain}},'' {\em IEEE Transactions on Automatic Control},
  vol.~67, no.~10, pp.~5570--5577, 2022.

\bibitem{Selivanov2023b}
A.~Selivanov and E.~Fridman, ``Finite-{{Dimensional Boundary Control}} of a
  {{Wave Equation With Viscous Friction}} and {{Boundary Measurements}},'' {\em
  IEEE Transactions on Automatic Control}, pp.~1--8, 2023.

\bibitem{Selivanov2023}
A.~Selivanov and E.~Fridman, ``Disturbance {{Attenuation}} in the
  {{Euler}}–{{Bernoulli Beam}} with {{Viscous}} and {{Kelvin}}–{{Voigt
  Damping}} via {{Piezoelectric Actuators}},'' in {\em 62nd {{IEEE Conference}}
  on {{Decision}} and {{Control}}}, pp.~1961--1966, 2023.

\bibitem{Bontsema1988a}
J.~Bontsema, R.~F. Curtain, and J.~M. Schumacher, ``Robust control of flexible
  structures {{A}} case study,'' {\em Automatica}, vol.~24, no.~2,
  pp.~177--186, 1988.

\bibitem{Lenz1993}
K.~Lenz and H.~Özbay, ``Analysis and {{Robust Control Techniques}} for an
  {{Ideal Flexible Beam}},'' in {\em Control and {{Dynamic Systems}}} (C.~T.
  Leondes, ed.), vol.~57 of {\em Multidisciplinary {{Engineering Systems}}:
  {{Design}} and {{Optimization Techniques}} and Their {{Application}}},
  pp.~369--421, Academic Press, 1993.

\bibitem{Halim2001}
D.~Halim and S.~O. Moheimani, ``Spatial resonant control of flexible structures
  - application to a piezoelectric laminate beam,'' {\em IEEE Transactions on
  Control Systems Technology}, vol.~9, no.~1, pp.~37--53, 2001.

\bibitem{Tucsnak1996}
M.~Tucsnak, ``Regularity and {{Exact Controllability}} for a {{Beam}} with
  {{Piezoelectric Actuator}},'' {\em SIAM Journal on Control and Optimization},
  vol.~34, no.~3, pp.~922--930, 1996.

\bibitem{Crepeau2006}
E.~Crépeau and C.~Prieur, ``Control of a clamped-free beam by a piezoelectric
  actuator,'' {\em ESAIM - Control, Optimisation and Calculus of Variations},
  vol.~12, no.~3, pp.~545--563, 2006.

\bibitem{Bai2024}
Y.~Bai, C.~Prieur, and Z.~Wang, ``Exact controllability for a {{Rayleigh}} beam
  with piezoelectric actuator,'' {\em Systems \& Control Letters}, vol.~186,
  p.~105759, 2024.

\bibitem{Halim2002}
D.~Halim and S.~Moheimani, ``Experimental implementation of spatial
  {{H-infinity}} control on a piezoelectric-laminate beam,'' {\em IEEE/ASME
  Transactions on Mechatronics}, vol.~7, no.~3, pp.~346--356, 2002.

\bibitem{Belyaev2018}
A.~Belyaev, A.~Fedotov, H.~Irschik, M.~Nader, V.~Polyanskiy, and N.~Smirnova,
  ``Experimental study of local and modal approaches to active vibration
  control of elastic systems,'' {\em Structural Control and Health Monitoring},
  vol.~25, no.~2, 2018.

\bibitem{Russell1992}
D.~L. Russell, ``On {{Mathematical Models}} for the {{Elastic Beam}} with
  {{Frequency-Proportional Damping}},'' in {\em Control and {{Estimation}} in
  {{Distributed Parameter Systems}}}, pp.~125--169, {Society for Industrial and
  Applied Mathematics}, 1992.

\bibitem{Herrmann2008}
L.~Herrmann, ``Vibration of the {{Euler-Bernoulli Beam}} with {{Allowance}} for
  {{Dampings}},'' in {\em World {{Congress}} on {{Engineering}}}, pp.~901--904,
  2008.

\bibitem{Crawley1987}
E.~F. Crawley and J.~De~Luis, ``Use of piezoelectric actuators as elements of
  intelligent structures,'' {\em AIAA Journal}, vol.~25, no.~10,
  pp.~1373--1385, 1987.

\bibitem{Fuller1996}
C.~R. Fuller, S.~J. Elliott, and P.~A. Nelson, {\em Active {{Control}} of
  {{Vibration}}}.
\newblock Elsevier, 1996.

\bibitem{Chen1982}
G.~Chen and D.~L. Russell, ``A mathematical model for linear elastic systems
  with structural damping,'' {\em Quarterly of Applied Mathematics}, vol.~39,
  no.~4, pp.~433--454, 1982.

\bibitem{Pazy1983}
A.~Pazy, {\em Semigroups of {{Linear Operators}} and {{Applications}} to
  {{Partial Differential Equations}}}. 
\newblock Springer, 1983.

\bibitem{Tucsnak2009}
M.~Tucsnak and G.~Weiss, {\em Observation and {{Control}} for {{Operator
  Semigroups}}}.
\newblock Birkhäuser, 2009.

\bibitem{Timoshenko1955}
S.~Timoshenko, {\em Strength of {{Materials}}. {{Part}} 1}.
\newblock Van Nostrand, 1955.

\bibitem{Banks1991}
H.~T. Banks and D.~J. Inman, ``On {{Damping Mechanisms}} in {{Beams}},'' {\em
  Journal of Applied Mechanics}, vol.~58, no.~3, pp.~716--723, 1991.

\end{thebibliography}
\end{document}